\documentclass{amsart}

\usepackage{amssymb}
\usepackage{graphicx}
\usepackage{hyperref}

\usepackage[margin=1in]{geometry}  
\usepackage{graphicx}              
\usepackage{amsmath}               
\usepackage{amsfonts}              
\usepackage{amsthm}                

\begin{document}
\pagenumbering{arabic}
\title{Exploration of another Sol Lewitt puzzle from Barry Cipra}

\author{Barry Cipra}
\author{Donna Dietz}
\author{Peter Winkler}

\begin{abstract}
  At MOVES 2019, Barry Cipra casually introduced a new
  ``Sol Lewitt'' puzzle to fellow conference goers.  Several
  brainstorming sessions ensued with Barry, Peter Winkler , Donna
  Dietz, and other attendees.  This paper is to document the puzzle
  and some insights so others can enjoy and build on this lovely
  puzzle.  (Look for it in an upcoming book by Peter!)
\end{abstract}

\maketitle

\section {Introduction}
Barry Cipra\footnote{bcipra@rconnect.com} has long been inspired by
the artwork of Sol Lewitt, particularly his ``Fifteen Etchings'' work
from 1973. In this piece of work a square is subdivided into 16
smaller squares, each with an etching (except for one blank
spot). These etchings contain horizontal, vertical, and diagonal
lines. In Barry's first Sol Lewitt puzzle, the goal was to rearrange
the smaller squares, so that no lines terminated inside the 4x4
grid. This puzzle and its solutions can be explored online at many
sites including \cite{website1}.

\begin{figure}
\begin{center}
\includegraphics[width=3in]{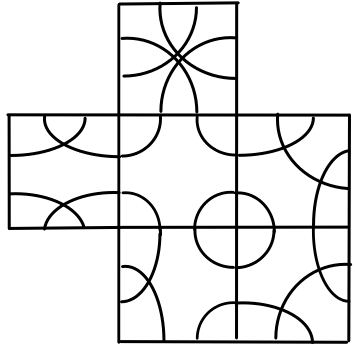}\\
\end{center}
\caption{No paths terminate inside the tiled region.}
\label{demo1}
\end{figure}

As a takeoff on this original puzzle, Barry designed 16 new square
puzzle tiles filled with curved segments, such that Cartesian
placements of the tiles will never result in curves termanating within
the tiled region as demonstrated in Figure \ref{demo1}. If you follow
any curve in either direction, it will either be part of a loop, or
both ends will termanate at the border.

\begin{figure}
\begin{center}
\includegraphics[width=5in]{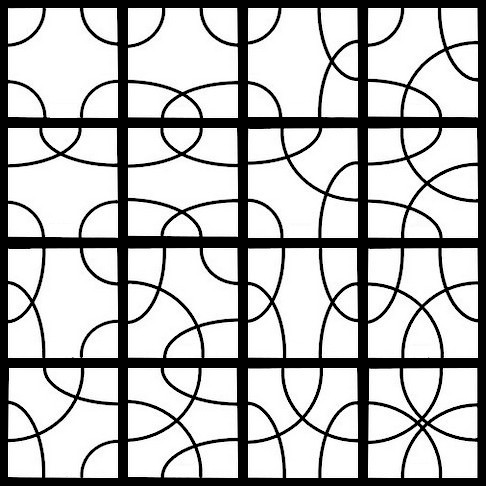}\\
\end{center}
\caption{Puzzle tiles.}
\label{UnnumberedTiles}
\end{figure}

Barry brought these tiles to MOVES 2019, asking fellow attendees
(including Peter Winkler\footnote{peter.winkler@dartmouth.edu} and
Donna Dietz\footnote{dietz@american.edu}) to ponder them and tell him
what they thought. In the true spirit of mathematical adventure, we
encourage readers to immediately stop reading, and go play with these
tiles either on paper by cutting out Figure \ref{UnnumberedTiles} with
scissors or online at \cite{website2}, before returning to read about
what we discovered.

\section{The design}
The tiles themselves can be thought of as rooms with four entrances,
one each side of the tile, allowing access to the
center of the room.  Now also imagine curtains for each entrance which
could be left open or closed. When open, they sit on either side of
the entrance, but when closed will overlap. Another way to visualize
this is to see that each side of each tile involves either just the
two points on the side, or it involves a third point, where the curves
cross.  On each side, we choose whether to activate that crossing
point or not.  Since there are 16 ways to do this, we end up with 16
oriented tiles.

\begin{figure}
\begin{center}
\includegraphics[width=3in]{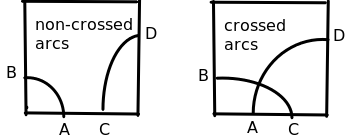}
\end{center}
\caption{Crossed and Non-Crossed Arcs}
\label{CrossedNonCrossedArcs}
\end{figure}

We can also focus on the curves rather than the spaces made by the curves.
Each curve is made up of arcs which deflect 90 degrees as they pass
though the tile.  Although there are different physical lengths, we
think of each arc as having length one.  There are binary strings or
numbers associated with each tile. The digits represent the crossing
or non-crossing (Figure \ref{CrossedNonCrossedArcs}) nature of each
side of the tile, starting with the left side of the tile and working
clockwise.

\section{A loopy question}

After a bit of exploration, we agreed to place the tiles on a 4x4 grid
representing a torus, and to preserve the original orientations of the
tiles. (This is by no means necessary, and readers are encouraged to
explore as many variations to this puzzle as they can create!)

When placed on a 4x4 torus, the far right and left edges are one
single edge, just as the top and bottom edge are the same edge.  (This
is the same structure as a Pac-Man board or a standard doughnut has.)
There are now no edges so curves cannot terminate.  The only option
left for any curve is to be a loop!

The loops may be strictly planar, or they may circle around the torus
once or more either horizontally or vertically.

\begin{figure}
\begin{center}
\includegraphics[width=3in]{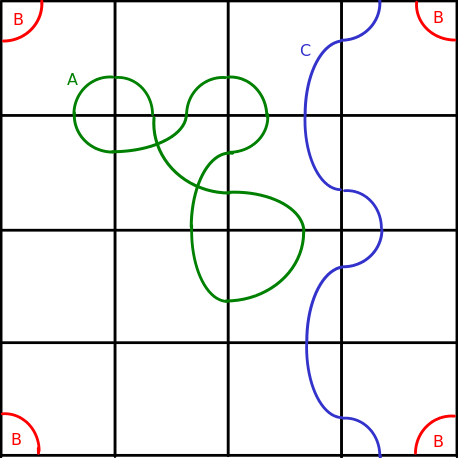}\\
\end{center}
\caption{Three loops.}
\label{demo2}
\end{figure}

In Figure \ref{demo2}, Curve A could just as easily have been created
on a plane as on a torus, so we will call it a planar loop.  Curve B
is also a planar loop, and it can be drawn as a circle when the torus
is cut differently.  Curve C, however, wraps around the torus in order
to form the loop. We will refer to these as planar loops and torus
loops respectively.

\begin{figure}
\begin{center}
\includegraphics[width=3in]{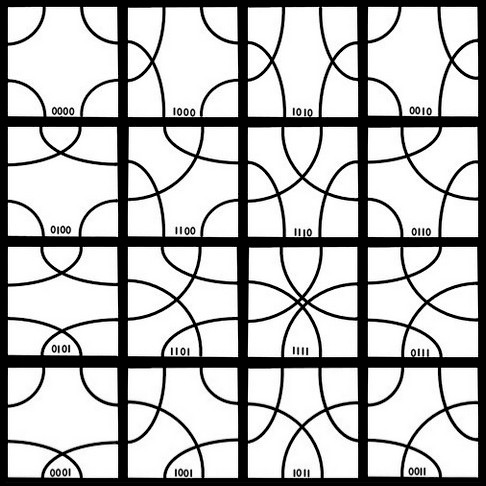}\\
\end{center}
\caption{Eight torus loops}
\label{BarrysSnakes}
\end{figure}

\begin{figure}
\begin{center}
\includegraphics[width=3in]{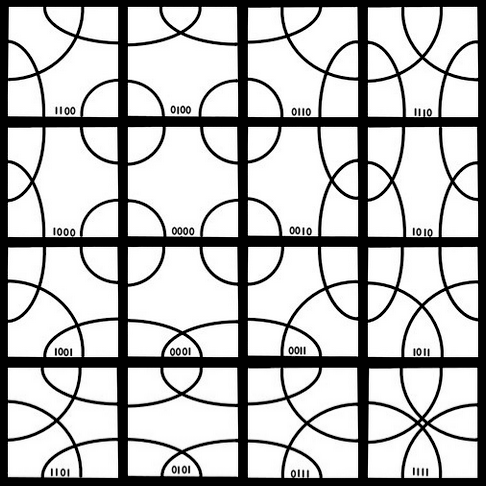}\\
\end{center}
\caption{Sixteen planar loops}
\label{loops}
\end{figure}

The first questions which arose between us concerned the maximal and
minimal numbers of loops possible and their lengths.  Experience had
led to several hunches, and a few configurations stood out as
particularly interesting. Figure \ref{BarrysSnakes} demonstrates how
it is possible to have 8 torus loops on the torus, while Figure
\ref{loops} demonstrates how it is possible to have 16 planar loops.

Given our numbering convention, these two configurations have
something special in common.  We can find them by using Gray
Code. Gray Code is named after Frank Gray, and is discussed in several
places online, such as \cite{FG}. Gray Code allows binary numbers to
be listed one after another and returning to the starting point, from
0 to $2^n -1$ while only changing the value of one digit at a
time. This efficiency can pay dividends in many applications.  Here,
it creates aesthetically pleasing patterns.  Both Figure
\ref{BarrysSnakes} and Figure \ref{loops} exhibit this property in
both the horizontal and vertical directions. Other Gray Code
arrangements exist for these tiles.  The reader is encouraged to
explore and enjoy them.

The configuartion in Figure \ref{BarrysSnakes} shows how the maximum number
of torus loops may be attained.  Since there are 64 arcs in total, and
we note that it takes 8 arcs to reach around the torus in either
direction, there cannot be more than 8 torus loops.

The configuration in Figure \ref{loops} shows how the maximum number of loops
(all plane loops) can be attained.  Since arcs deflect 90 degrees, you
cannot make any loop with fewer than four arcs, thus there cannot be
more than 16 loops.

\begin{figure}
\begin{center}
\includegraphics[width=3in]{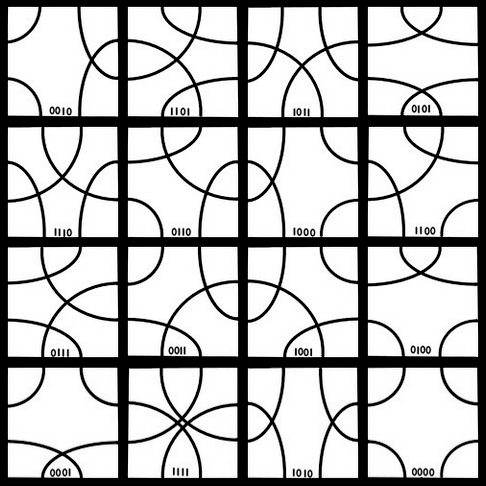}\\
\end{center}
\caption{A loop of length 60 and one of length 4}
\label{sixty}
\end{figure}

These configurations, as well as the one in Figure \ref{sixty}
were found by manual experimentation with the puzzle pieces.  Figure
\ref{sixty} demonstrates a loop of length 60 together with a loop
of length four.

\section{Parity of number of loops matches parity of crossings}

The question of course was whether or not a single loop of length 64
could be formed.  We suspected this was not possible, because all
configurations of the tiles had resulted in an even number of loops!
After several attempts at a proof, a simple, inductive case-by-case
proof was found by Peter.

\newtheorem{theorems}{Theorem}

\begin{theorems}
(Parity Theorem)
The number of loops and the number of crossings have the same parity.
\end{theorems}  

As a base case, take any configuration. However,  since it is not necessary to restrict the use of tiles, tile 0000 may be used in all 16 positions so there are initially no crossings in the configuration, making it easy to visualize.

\begin{figure}
\begin{center}
\includegraphics[width=3in]{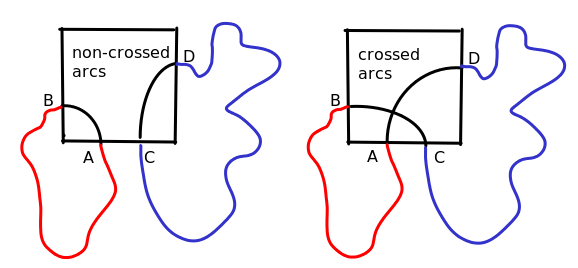}\\
\end{center}
\caption{Case A}
\label{CaseA}
\end{figure}

Case A) As in Figure \ref{CaseA}, there could be a tile with
non-crossed arcs coming through adjacent corners, which (perhaps
unbeknownst to us) are each part of two separate loops. If we replace
that set of arcs with crossed arcs using the same attachment points on
the tiles, this will result in a single loop.  We have decreased the
loops by one and increased the crossings by one, thereby maintaining
parity.

\begin{figure}
\begin{center}
\includegraphics[width=3in]{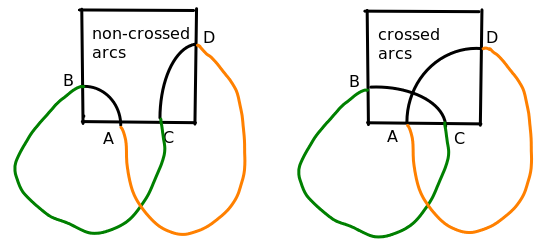}\\
\end{center}
\caption{Case B}
\label{CaseB}
\end{figure}

Case B) As in Figure \ref{CaseB}, there could be a tile which starts
out as a single loop, making a figure eight crossing outside of the
tile, but after the switch it becomes two loops. This will become two
separate loops as shown. We have increased the number of crossings by
one and also increased the number of loops by one thereby maintaing
parity.

Case C) This case is special, because it is a red herring. We will
show that this case does not exist for this set of tiles, given the
rules we have agreed upon. However it would appear that if there is a
connection between the lower to attachment points which does not
intersect the connection between the upper attachment points, when the
non-crossed arcs are replaced with crossed arcs, we should increase
only the crossings without changing the number of loops.
In order to explain why Case C cannot happen, we need another theorem.

\begin{theorems}
  \label{foo}
  (Checkerboard Crossing Theorem)
Once the direction is set for a loop, this induces a checkerboard pattern on all the tiles, constraining only horizontal entry (H) into the tile or only vertical entry into the tile (V).
\end{theorems}

\begin{figure}
\begin{center}
\includegraphics[width=3in]{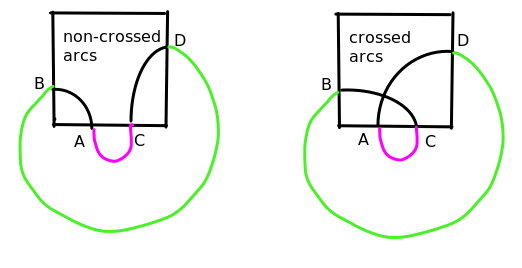}\\
\end{center}
\caption{Case C}
\label{CaseC}
\end{figure}

\begin{figure}
\begin{center}
\includegraphics[width=3in]{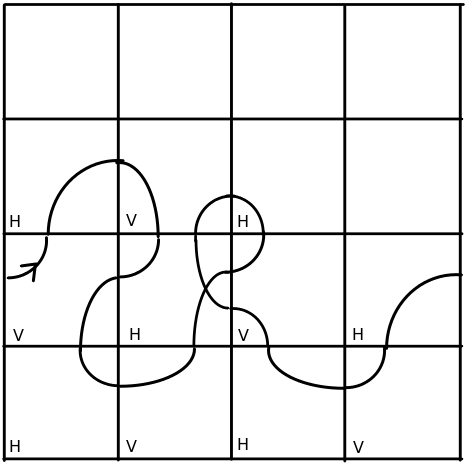}\\
\end{center}
\caption{Checkerboard tiling example}
\label{HVdemo}
\end{figure}

\begin{figure}
\begin{center}
\includegraphics[width=1in]{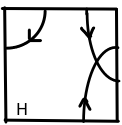}\\
\end{center}
\caption{Zoom-in on one tile}
\label{HVdemo2}
\end{figure}

When the loop enters a tile, that arc must come from either a
horizontal or a vertical side.  Since all arcs turn a quarter-circle,
it must exit opposite to how it entered thereby entering the next tile the other
way.  Our torus is 4x4, or 8x8 if you were to count arc steps.
Therefore, the entire torus (everywhere the curve may wander) has to
honor this checkerboard rule.  In Case C (Figure \ref{CaseC}), in
order for this case to happen, the directions of the arcs must be
opposite, which is not possible.  A loop is given in Figure
\ref{HVdemo} to help the reader visualize why this is true. A zoom-in
of one tile is given in Figure \ref{HVdemo2} indicating the local
directional structure induced on that tile as a result of choosing the
curve direction.

So, although the curve direction is arbitrary, and therefore the H or
V designation is arbitrary, as soon as a direction is chosen for a
curve, the whole curve must abide by that directional choice.

In our original tile portfolio, we had one of each kind of tile, which
means we had an even parity on the curve crossings.  Therefore we have
to maintain an even number of loops, making 60 the largest possible
loop, as demonstrated in Figure \ref{sixty}.

\section{Possible loop lengths}

It may appear obvious that the loops must have a length divisible by
four, but precisely because it seems obvious, it would
seem prudent to prove this fact.  We have two proofs, in fact!

\begin{figure}
\begin{center}
\includegraphics[width=3in]{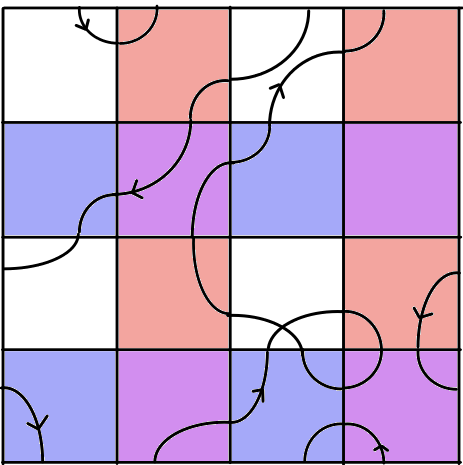}\\
\end{center}
\caption{Path length divisibility proof using four colors}
\label{CurvyCurveColored}
\end{figure}

The first proof is is short and sweet, and it makes use of the fact that
we have agreed to use a torus made of even numbers of tiles in each
direction.  Since we alternate our tile crossings as
vertical-horizontal-vertical-horizontal, we must alternate switching the
column, then row, then column, then row. This means we have four kinds
of tiles which must be visited in a specific order. In Figure
\ref{CurvyCurveColored}, this color sequence is {\it white, red,
  purple, blue}.  Since we visit all colors equally, the total path length is divisible by four.

The second proof is not very elegant, but it is true without the
restriction that the torus must have even numbers of tiles in
each direction.

Since arcs turn through quarter circles, it is obvious that the
minimal number of arcs in a loop must be four. Also obvious, due to
argments in the proof above (Theorem \ref{foo}), the loop must have an even
length. One possible attack would be to take existing loops and join
them as follows.  Take the entire torus and cover it with tile 0000,
thus making 16 loops.

\begin{figure}
\begin{center}
\includegraphics[width=4in]{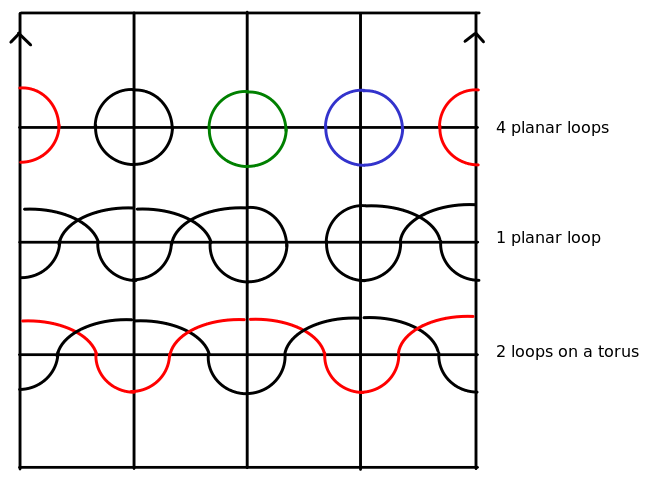}\\
\end{center}
\caption{Changing four original loops}
\label{planarTorusLoops}
\end{figure}

One by one, replace tiles with tiles containing just one more arc
crossing, until the desired configuration is attained. One example of
such a transition is shown in Figure \ref{planarTorusLoops}.  However,
even though this makes sense intuitively, the actions on the loops as
they split and join are quite subtle, so we will focus on the analysis
of an existing loop rather than creating loops.

All loops can be decomposed into their components based on their
topology. Each type of loop must satisfy the property that the length
is a multiple of four. When joining or splitting loops, we have seen
that the number of edges remains a constant.

Type I) A planar loop with some integer number of full circular turns
will have a length divisible by four, provided any net-zero deflection
loop also does.  This is because each full turn will require four arcs.

Type II) A torus loop in either direction will also have a length
divisible by four, provided any net-zero deflection loop also
does. This is because the arc-height and arc-width of our torus is 8,
which is divisible by four.

Type III) Any net-zero deflection loop in the plane must have length
divisible by four. The simplest case for this type of loop is a figure
eight. (This will be demonstrated in Theorem \ref{foo2}.)

Type IV) Any other loops must be decomposible into the previous types.

\begin{theorems}
  \label{foo2}
(Net zero deflection Theorem)
Zero net-deflection loops in a plane have length divisible by four.
\end{theorems}

\begin{figure}
\begin{center}
\includegraphics[width=4in]{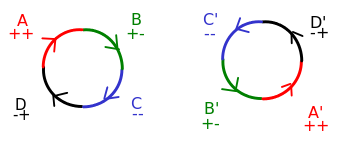}\\
\end{center}
\caption{Eight arc types}
\label{ABCD}
\end{figure}

There are four types of arcs we draw on the page, but each can have
two directions.  Thus, we will discuss all eight arcs as if they are
different.  Figure \ref{ABCD} gives the eight arc types, along with +
and - symbols indicating their x and y displacements.

The symbols A, B, C, D, A', B', C', and D' will be used not only to
represent the eight arc types, but also to indicate the number of times
that arc appears in a loop.  Thus, arc A appears A times in our loop.

\begin{figure}
\begin{center}
\includegraphics[width=2in]{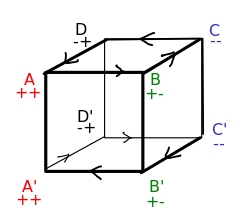}\\
\end{center}
\caption{Permitted transitions}
\label{Transitions}
\end{figure}

Due to the rules we've set for this game, not all arc types may follow
each other. In fact, once an arc has been chosen, only two
possibilites follow.  If the arc is at a point of inflection, it will
continue in the same direction but the concavity will change.
Otherwise, it must move to the next direction permitted. There must be
an even number of points of inflection in a loop. The permitted
transitions are shown in Figure \ref{Transitions}.  Note that this is
a cube and is therefore bipartite. Eight edges are directional, while
four edges can be traversed in either direction.  The cube's corners
represent arcs in the loop, while the edges represent points in the
loop.

Since the loop must follow this bipartite graph, and since we know from Theorem \ref{foo} that there must be an even number of arcs in the loop, we have this {\it Bipartite loop condition} as follows:

$$A+C+B'+D' = A'+C'+B+D. $$

Since we have a net zero deflection loop, we have a {\it Zero deflection condition}:

$$A + B+ C + D = A'+B'+C'+D'.$$

Finally, since we have to balance our + and - tallies on both x and y coordinates, we have two more conditions, the {\it Balancing conditions}:

$$ A +A'+ B+B' = C+C'+D+D', and$$
$$A+A'+D+D' = B+B'+C+C'.$$

Working with just the balancing conditions, subtracing them from each other gives the {\it Simplified balancing conditions}:

$$A+A' = C+C', and$$
$$B+B' = D+D'.$$

If the Bipartite loop condition is subtracted from the zero deflection condition (flipped), we get that $$(A+C) - (A'+C') = (A'+C') - (A+C),$$  which tells us that

$$A+C = A'+C', and\; similarly$$
$$B+D = B'+D'.$$

But this reveals something interesting.  Let $A+A'=C+C'=2x$, $A+C=A'+C'=2y$ which means that $x=y$. This gives the {\it Solved Balancing Equations}:

$$A=A'=C=C', and$$
$$B=B'=D=D'.$$

\begin{figure}
\begin{center}
\includegraphics[width=4in]{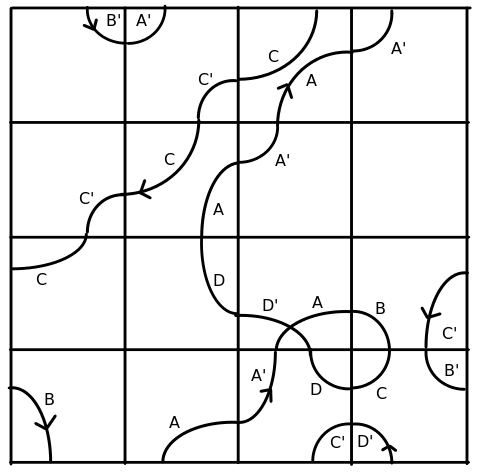}\\
\end{center}
\caption{Example curve labeled by arc types}
\label{CurvyCurve}
\end{figure}

Since the total number of arcs in a zero deflection plane curve is $A+B+C+D+A'+B'+C'+D'$, or $4x+4y$, it is divisible by four! Figure \ref{CurvyCurve} shows how this works out for one specific loop, where $A=A'=C=C'=4$ and $B=B'=D=D'=2$.

\section{Right of way}

\begin{figure}
\begin{center}
\includegraphics[width=1in]{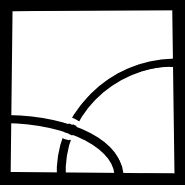}\\
\end{center}
\caption{Right of way determination for two arcs}
\label{RightOfWay1}
\end{figure}

\begin{figure}
\begin{center}
\includegraphics[width=1in]{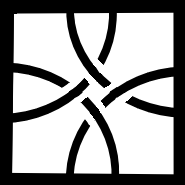}\\
\end{center}
\caption{Four arc crossings in one tile}
\label{RightOfWay2}
\end{figure}

As we explored this set of game tiles, it occured to us that it might be nice to indicate a three dimensionality on the tiles, so one arc would appear to be above and one below at each crossing.  We decided to use {\it right of way} to determine which arc would get the upper level.  So, as shown in Figure \ref{RightOfWay1}, as two arcs come in from the bottom of the tile, the one on the right has priority. A maximum of four crossings may occur on one tile, such as shown in Figure \ref{RightOfWay2}.

\begin{figure}
\begin{center}
\includegraphics[width=5in]{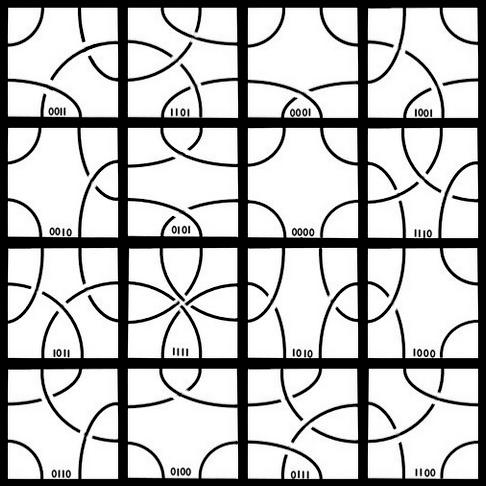}\\
\end{center}
\caption{An entire torus marked with right of way bridges}
\label{snakesOnAPlane}
\end{figure}

An example of the entire torus being marked this way is shown in Figure \ref{snakesOnAPlane}.

You may notice in Figure \ref{snakesOnAPlane}, that there is a wonderful weaving motion, over-under-over-under, on any of these loops.  Why is that?

\begin{theorems}
(Weave Theorem)
  Tiles marked with right-of-way bridges will always alternate over-under-over-under.
\end{theorems}

\begin{figure}
\begin{center}
\includegraphics[width=1.5in]{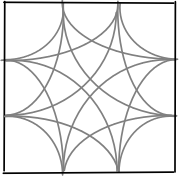}\\
\end{center}
\caption{All possible arcs drawn together }
\label{masterStar}
\end{figure}

\begin{figure}
\begin{center}
\includegraphics[width=1.5in]{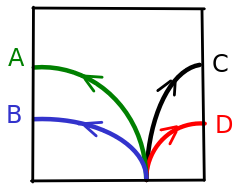}\\
\end{center}
\caption{All possible arcs starting at one entry point}
\label{preStar}
\end{figure}

First, let us consider all possible crossings, by simply drawing all possible arcs on one single tile as in Figure \ref{masterStar}, and let us now focus on what happens only upon entry up from the right side, on a horizontal side as in Figure \ref{preStar}.

\begin{figure}
\begin{center}
\includegraphics[width=1.5in]{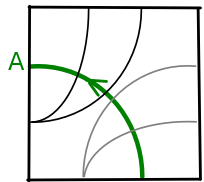}\\
\end{center}
\caption{Intersections when exiting at A}
\label{Astar}
\end{figure}

In the case of exiting at point A, consider the gray curves and the black curves from Figure \ref{Astar}.  We will have one crossing from each set of curves, thus we'll enter the tile with right of way, cross above, then below, then retain right of way as we exit, entering the next tile also on the right.

In the cases of exiting at points B or C from the tile, we will enter the tile on the right with right of way, then use our right of way once, but leaving the tile on the left, and without right of way.

In the final case, exiting at point D, we enter and exit with right of way, like in A, but we've used zero crossings, which is fine, as zero is just as even as two is.

A similar analsis can be accomplished on the four remaining paths of
travel through the tile.  Thus, weaving will always be maintained as
we traverse our loops in this puzzle.

\begin{figure}
\begin{center}
\includegraphics[width=1.5in]{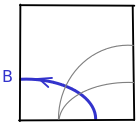}\\
\end{center}
\caption{Intersections when exiting at B}
\label{Bstar}
\end{figure}

\begin{figure}
\begin{center}
\includegraphics[width=1.5in]{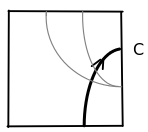}\\
\end{center}
\caption{Intersections when exiting at C}
\label{Cstar}
\end{figure}

\begin{figure}
\begin{center}
\includegraphics[width=1.5in]{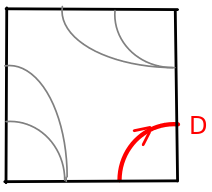}\\
\end{center}
\caption{Intersections when exiting at D}
\label{Dstar}
\end{figure}

\section{Edge loops}
\begin{figure}
\begin{center}
\includegraphics[width=3in]{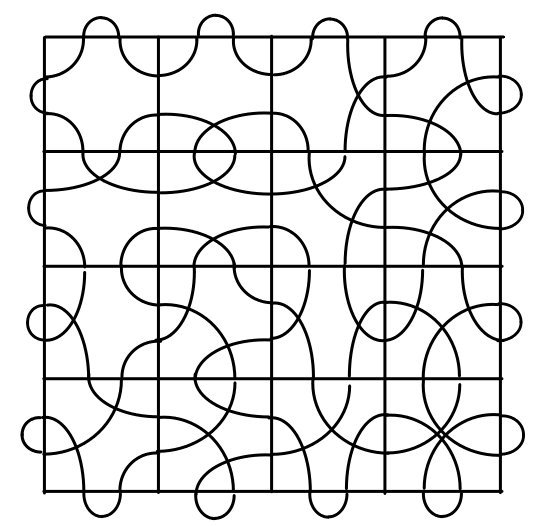}\\
\end{center}
\caption{Another variation}
\label{JimPropp}
\end{figure}

One variation idea by Jim Propp was to not use a torus, but to put loops on all the edge pieces as shown in Figure \ref{JimPropp}. This certainly changes the puzzle in a dramatic way!

\section{Magic Squares}
\begin{figure}
\begin{center}
\includegraphics[width=3in]{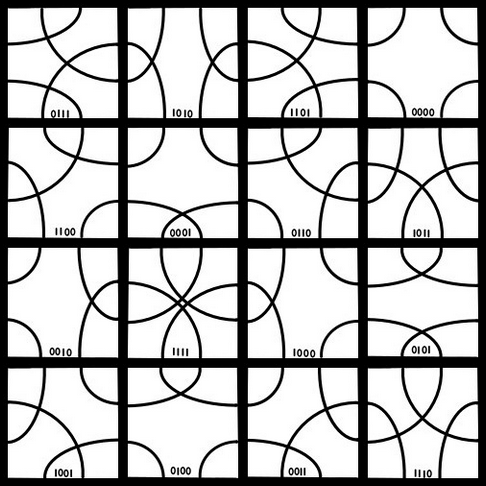}
\end{center}
\caption{A magical variation}
\label{magic}
\end{figure}

Another variation is to consider the binary labels as binary numbers, and then look at arrangements that form pandiagonal magic squares. Figure \ref{magic} shows one such arrangement.

\section {Conclusions}
In conclusion, we feel this is a fun set of tiles, and worthy of more exploration, more questions, and more theorems. Perhaps we may add more restrictions or loosen them.  Feel free to explore and touch base with the authors about your findings!

\end{document}